\documentclass[paper,onecolumn,a4paper]{IEEEtran}

\usepackage{cite}
\usepackage{amsmath,amssymb,amsbsy}
\newtheorem{theorem}{Theorem}[section]
\newtheorem{lemma}{Lemma}[section]
\newtheorem{definition}{Definition}[section]

\newtheorem{corollary}{Corollary}[section]

\title{Optimal H2 order-one reduction by solving eigenproblems for
polynomial equations
}

\author{
Bernard~Hanzon, Jan~M.~Maciejowski
and~Chun Tung Chou%
\thanks{Corresponding author Prof. J.M. Maciejowski.
Tel +44 1223 332732 Fax +44 1223 332662; E-mail jmm@eng.cam.ac.uk.}%
\thanks{B. Hanzon is with the School of Mathematical Sciences,
University College Cork, Cork City, Ireland.}%
\thanks{J. Maciejowski is with the Cambridge University Engineering
  Department, Trumpington Street, CB2 1PZ Cambridge, England.}%
\thanks{C.T. Chou is with the School of Computer Science and Engineering,
The University of New South Wales, Sydney, NSW2052, Australia.}%
}

\begin{document}

\maketitle

\begin{abstract}
A method is given for solving an optimal $H_2$ approximation
problem for SISO linear time-invariant stable systems. The method,
based on constructive algebra, guarantees that the global optimum
is found; it does not involve any gradient-based search, and hence
avoids the usual problems of local minima. We examine mostly the
case when the model order is reduced by one, and when the original
system has distinct poles. This case exhibits special structure
which allows us to provide a complete solution.  The problem is
converted into linear algebra by exhibiting a finite-dimensional
basis for a certain space, and can then be solved by eigenvalue
calculations, following the methods developed by Stetter and
M\"{o}ller~\cite{MollerStetter,Stetter}. The use of Buchberger's
algorithm is avoided by writing the first-order optimality
conditions in a special form, from which a Gr\"{o}bner basis is
immediately available. Compared with our previous work
~\cite{HanMacAut96}, the method presented here has much smaller
time and memory requirements, and can therefore be applied to
systems of significantly higher McMillan degree. In addition, some
hypotheses which were required in the previous work have been
removed. Some examples are included.
\end{abstract}

\section{Introduction}

In this paper we consider the problem of approximating a stable
linear dynamic system by one of lower McMillan degree. We take the
$L_2$ norm as the measure of approximation, namely we solve the
problem
\begin{equation}
\min_{\hat{h}\in {\cal M}(n)} \int_0^\infty |h(t)-\hat{h}(t)|^2 dt
\end{equation}
where $h\in {\cal M}(N)$ is the impulse response of the original system,
$\hat{h}$ is the impulse response of the approximating system, and ${\cal
M}(N)$
denotes the set of impulse responses of minimal stable systems of
McMillan degree $N$. This problem is equivalent to the problem of finding an
approximation which minimizes the $H_2$ norm of the error in the frequency
response:
\begin{equation}
\min_{\hat{H}\in {\cal H}(n)} \frac{1}{2\pi}
\int_{-\infty}^\infty
|H(\omega)-\hat{H}(\omega)|^2 d\omega
\end{equation}
where $H$ and $\hat{H}$ are the frequency responses of the
original and the approximating systems, respectively, and ${\cal
H}(N)$ denotes the set of Fourier transforms of elements of ${\cal
M}(N)$. Throughout this paper we consider SISO systems only, and
we solve the $H_2$ problem for $n=N-1$. We assume mostly that the
`true' system has distinct poles.
From section II onwards we will
work with the set $\Sigma S_N$ of rational transfer functions,
whose impulse responses are elements of ${\cal M}(N)$ and
frequency responses are elements of ${\cal H}(N)$, and we will
look for approximants in the set $\Sigma S_n$.

The $H_2$ problem has many applications and connections to other
problems in systems and control theory, including model
simplification, system identification, and approximate model
matching. Many publications treat this problem, such as
\cite{MeiLuen,baratchartolivi88} and the references cited therein.
An early publication on this problem, possibly the oldest, is
\cite{aigrainwilliams}. We investigate the $H_2$ approximation
problem by means of constructive algebra, in particular by
exploiting the theory of polynomial ideals. There is an increasing
use of computer algebra in systems theory, see e.g.\
\cite{FliessGlad,Helton,Oberst,RochaWillems,WalGeoTan01,WoodRogersOwens}.
This paper makes a further contribution to this trend.

We believe that the significance of this paper lies in its
introduction of a promising new approach to model reduction
problems. We emphasise that this approach does not involve
gradient-based search methods, and hence avoids the usual problems
associated with local minima.  Our use of constructive algebra
leads to an algorithm with the important attribute that the
solution found is guaranteed to be the global optimum. In
\cite{HanMacAut96} two of the present authors already applied
constructive algebra to the $H_2$ approximation problem, taking an
approach based on state-space realizations of the linear systems
involved. By contrast, the approach here is based on a form of the
first-order necessary conditions for optimality which arises from
transfer function descriptions of both the original and the
approximating systems. The solution method which we develop here
is quite different from that developed in \cite{HanMacAut96}.
Computationally it is much more efficient, as regards both memory
and time requirements.  This allows us to tackle problems with
significantly larger values of $N$, as can be seen from the
examples. Furthermore, \cite{HanMacAut96} required some technical
hypotheses relating to the finiteness of the number of critical
points, which are not needed in this paper.

In addition to finding the global optimum, our approach gives
important new theoretical insight into the structure of the
reduction-by-one problem. In particular, we show that the number
of critical points is finite, and in fact no greater than $2^N-1$.
The computational complexity is high, and the method involves some
delicate numerical steps, so we do not claim that our approach is
a rival, at this stage, for conventional numerical approaches in
routine applications to engineering problems. But even now it has
some practical uses, for example as a generator of `benchmark'
solutions against which other methods can be tested. Since, as
will be seen, it relies on eigenvalue calculations for a set of
matrices which can be constructed in a rather straightforward
manner, our approach is in some ways comparable with Glover's
method for solving the Hankel-norm approximation problem
\cite{Glover84}. Promising developments which combine the current
approach with numerical methods for solving large eigenvalue
problems in related applications are reported in
\cite{BleHanPee03}.

In the next section we obtain a special representation of the first-order
necessary conditions for optimality. This representation is in the form of
a set of quadratic equations, which take a special form which we call
{\em diagonal quadratic}.
The following section investigates such diagonal quadratic equations.
It is shown that the polynomials which define these equations form a
Gr\"{o}bner basis for the ideal generated by themselves. It is further shown
that these equations have a finite set of solutions, and that in
consequence a certain space is finite-dimensional. Furthermore a basis for
this space is identified, which allows a solution method based on linear
algebra.

We then present such a method of solving a system of polynomial
equations.  This method relies on obtaining a Gr\"{o}bner basis,
but in the application to the specific $H_2$ problem considered
here, such a basis is immediately available. This method of
solving polynomial equations is of general use and it is known in
the computer algebra community, see
\cite{Stetter,MollerStetter,Corless} and the references therein.
The development here is self-contained and starts with
constructing a matrix solution of the system of polynomial
equations, from which the desired solutions can be found by
solving a collection of eigenvalue-eigenvector problems.
These eigenproblems can be solved either by numerical methods or by
symbolic methods.
We believe from a system theoretic
point of view it is very natural to start with the construction of a
matrix solution; in fact the matrices obtained are
generalised companion matrices.

A section then applies this method to the solution of the $H_2$
problem, for the case $n=N-1$ and distinct poles of the original
system. How to treat repeated poles is outlined in a short
section. This is followed by two examples.

\section{A special representation of the first order conditions.}
\label{foc}

In this section the first order conditions for a class of  $H_2$
model order reduction problems will be considered. Studying the
outcomes of a computer algebra calculation in which a set of
symbolic first order conditions for the $H_2$ model order
reduction problem was brought into a recursive form, it was
observed that the occurrence of multiple poles in the original
system gave rise to a certain singularity  in the first order
equations. This was the motivation for investigating the class of
systems with distinct poles separately from the class of systems
with multiple poles. The continuous-time case is treated here, but
the discrete-time case is in fact the same up to isometry (see
e.g. \cite{HanzonTracts}, Theorem 5.4-3; \cite{HanzonTUD}, Theorem
3.2-22).

Now let us set up the problem.
In fact there are several
equivalent formulations. One formulation which is closest
to the form of the first order conditions that we use
in this paper is as follows.(For other formulations refer to the
literature, e.g. \cite{HanMacAut96})

Consider a continuous-time stable SISO linear system.
Without loss of
generality we can assume the system to be strictly proper, because
 if it is not then the direct feedthrough term of the optimal
$H_2$ approximant will be equal to the direct feedthrough term
of the original system, and the strictly proper part of the optimal
 approximant will not be influenced at all (nor will the
 strictly proper part of any of
 the critical points) by the value of the
direct feedthrough term.
Let the transfer function of the  original system (i.e. the system
 that is to be reduced in order) be given
by $e(s)/d(s),$ where $e$ is some polynomial with real
 coefficients  of degree at most $N-1,$
and $d$ is a monic polynomial with real coefficients of degree $N$
 with all its zeroes (i.e. poles of the transfer function)
 $\delta_1,\delta_2,\ldots,\delta_N,$ within the open left half
of the complex plane. Assume that $e$ and $d$ are coprime.

Consider the rational function $\frac{e(s)}{d(s)}.$ It is an
element of the Hardy space $H_2$ of square summable functions on the
imaginary axis which are analytic on the open right halfplane and
satisfy a certain continuity requirement on the imaginary axis(cf.
\cite{Hoffman}).  In this paper we work with the subspace of real
rational functions in $H_2.$ This subspace consists of all strictly
proper real rational functions which have the property that all the
poles lie in the open left half plane. The space $H_2$ is in fact a
Hilbert space with corresponding norm $\|.\|_2$ of a function $t \in
H_2$ given by
\[\|t\|_2^2=\frac{1}{2 \pi}
\int_{-\infty}^{\infty} |t(i\omega)|^2 d\omega
\]

Consider the differentiable manifold $\Sigma S_{n}$ of all real
rational functions $\frac{b(s)}{a(s)}$ in $H_2$ such that $b(s)$ and
$a(s)$ are coprime, the coefficients of $a(s)$ and $b(s)$ are real and
$a(s)$ is a Hurwitz polynomial of degree $n.$ For more information
about the structure of this differentiable manifold see for example
\cite{ChouHanzonSCL} and \cite{HanzonOber97} and the references given
there.  The $H_2$ model order reduction problem
can now be formulated as the following optimization problem:
\[ \min_{\frac{b(s)}{a(s)} \in \Sigma S_{n}} \left\|
\frac{e(s)}{d(s)}-\frac{b(s)}{a(s)} \right\|_2.\]

{\em Remark.} It is well-known that the distance squared
$\left\| \frac{e(s)}{d(s)}-\frac{b(s)}{a(s)} \right\|_2^2$
is in fact a rational function of the coefficients of the
numerator and denominator polynomials (see the literature, e.g.
\cite{HanzonTracts}; in order to obtain explicit rational
 function formulas
one could use the methods proposed in \cite{HanzonPeetersLAA96} )

A well-known first order necessary condition for optimality of an
$n-$th order transfer function  $b(s)/a(s)$ with real
coefficients, as an approximant in $H_2$ is the following.
First let us present a geometric formulation.

If $\frac{b(s)}{a(s)}$ is an optimal approximant of the transfer
function  $\frac{e(s)}{d(s)}$ with respect to the  $H_2-$norm,
then  the difference $\frac{e(s)}{d(s)}-\frac{b(s)}{a(s)}$ is
perpendicular to the tangent plane at the manifold of
transfer functions of order $n$ at the point $\frac{b(s)}{a(s)}.$
\\
It is well-known (and not hard to show) that the tangent space
consists of all strictly proper rational functions of the form
$\frac{p(s)}{a(s)^2},$ where $p$ is a polynomial of degree at most
 $2n-1.$ From the theory of Hardy spaces it follows that the
orthogonal complement in $H_2$ of this vector space is given
by $a(-s)^2 H_2,$ i.e. all $H_2-$functions which can be
written as the product of the function $a(-s)^2$ and an arbitrary
$H_2$ function. Combining this with the first order conditions
given above, it follows that the numerator of the
difference $\frac{e(s)}{d(s)}-\frac{b(s)}{a(s)}$
has to be divisible by $a(-s)^2.$
(Cf. \cite{MeiLuen}, see also \cite{baratchartolivi88},
\cite{baratchartcardelliolivi}).
Algebraically this can be
written down as follows:

Let $n<N.$
If $\frac{b(s)}{a(s)}$ is an optimal approximant within the class
of transfer functions  of order $n$ in $H_2,$ of the transfer
function  $\frac{e(s)}{d(s)}$ in $H_2,$ with respect to the
$H_2-$norm,
then there exists a polynomial $q(s)$ of degree at most
$N-(n+1)$ such that
\begin{equation}
\label{focII}
e(s)a(s)-b(s)d(s)=a(-s)^2 q(s).
\end{equation}

Let us now specialise to the case in which $n=N-1$ and the
original system has distinct poles, i.e. the multiplicity of each
of the $N=n+1$ poles $\delta_1,\ldots,\delta_N$ is one. The rest
of this paper concentrates mostly on this case. Now the polynomial
$q(s)$ has degree zero, so it reduces to a constant $q(s)=q_0.$
The unknowns in the polynomial equation are the polynomials $b(s),
a(s)$ and the number $q_0.$ Although $q_0$ is only an auxiliary
variable we will not eliminate it. Note that once the polynomial
$a$ and the number $q_0$ are known, the polynomial $b$ follows
from the formula
\begin{equation}
\label{bequation}
b(s)=\frac{e(s)a(s)-q_0 a(-s)^2}{d(s)}.
\end{equation}
Substituting $s=\delta_i, i=1,\ldots,N$ in the polynomial
equation (\ref{focII}) one obtains:
\begin{equation}
\label{focIII}
e(\delta_i)a(\delta_i)=a(-\delta_i)^2 q_0, i=1,\ldots,N.
\end{equation}
Note that the polynomials appearing here do not depend on
the polynomial $b,$ due to the fact that $d(\delta_i)=0$
for each $i=1,\ldots,N.$
Further note that the possibility $q_0=0$ can be excluded on the grounds that
if $q_0=0$ then either $e(\delta_i)=0$ for some value of
$i\in \{1,\ldots,N\},$ which implies that there is pole-zero
cancellation in the original transfer function and the
order of the transfer function will be smaller than $N,$
which can be ruled out without loss of generality,
or otherwise it would follow that $a(s)=0$ in $N$ different
points, namely at $s=\delta_i,i=1,\ldots,N,$ which together
with the fact that $a$ has degree $n=N-1$ would imply that $a=0,$
which is
in contradiction with the assumption that $a$ is monic.
It follows that $q_0\not=0$ for each value of $q_0$ that corresponds to a
solution of the first order equations. Therefore 
 multiplying both sides of the
polynomial equation with $q_0$ the first order conditions can be rewritten as
\begin{equation}
\label{focIV}
e(\delta_i)a(\delta_i)q_0=\left(a(-\delta_i) q_0\right)^2,
 i=1,\ldots,N, q_0\not=0.
\end{equation}
The polynomial $a$ is monic, so $q_0$ is the leading coefficient
of the non-zero polynomial $\tilde{a}:=q_0 a.$ Using this notation the
first order equations can be rewritten as
\begin{equation}
\label{focV}
e(\delta_i)\tilde{a}(\delta_i)=\tilde{a}(-\delta_i)^2,
 i=1,\ldots,N,~\tilde{a}\not=0.
\end{equation}

The idea is now to consider this as an equation in the
unknowns $\tilde{a}(-\delta_i), i=1,\ldots,N.$ In order to
do this explicitly we need to express the sequence of numbers
$\tilde{a}(\delta_i), i=1,\ldots,N$ in terms of the sequence
of numbers
$\tilde{a}(-\delta_i), i=1,\ldots,N.$
This can be done by relating both sequences to
the coefficients $\tilde{a}_j, j=0,\ldots,N-1,$ of the polynomial
 $\tilde{a}(s)=\tilde{a}_{N-1} s^{N-1}+\tilde{a}_{N-2} s^{N-2}+
\ldots+
\tilde{a}_{0} s^0.$
Let $V(\delta_1,\ldots,\delta_N)$ denote the Vandermonde matrix
\begin{eqnarray}
V(\delta_1,\ldots,\delta_N)&:=&
\left(\begin{array}{ccccc}
1&\delta_1&\delta_1^2&\ldots&\delta_1^{N-1}\\
1&\delta_2&\delta_2^2&\ldots&\delta_2^{N-1}\\
\vdots&\vdots&\vdots& &\vdots\\
1&\delta_N&\delta_N^2&\ldots&\delta_N^{N-1}
\end{array}\right).
\end{eqnarray}

Using matrix-vector notation the following
linear relations are obtained:
\begin{eqnarray}
\left[\begin{array}{c}
\tilde{a}(\delta_1)\\
\vdots\\
\tilde{a}(\delta_{N})\end{array}\right]&=&V(\delta_1,\ldots,\delta_N)
\left[\begin{array}{c}
\tilde{a}_0\\
\vdots\\
\tilde{a}_{N-1}\end{array}\right]
\end{eqnarray}
and
\begin{eqnarray}
\label{atilantidel}
\left[\begin{array}{c}
\tilde{a}(-\delta_1)\\
\vdots\\
\tilde{a}(-\delta_{N})
\end{array}\right]&=&V(-\delta_1,\ldots,-\delta_N)
\left[\begin{array}{c}
\tilde{a}_0\\
\vdots\\
\tilde{a}_{N-1}
\end{array}
\right].
\end{eqnarray}
It follows that
\begin{eqnarray}
\left[\begin{array}{c}
\tilde{a}(\delta_1)\\
\vdots\\
\tilde{a}(\delta_{N})\end{array}\right]&=&V(\delta_1,\ldots,\delta_N)
V(-\delta_1,\ldots,-\delta_N)^{-1}
\left[\begin{array}{c}
\tilde{a}(-\delta_1)\\
\vdots\\
\tilde{a}(-\delta_N)
\end{array}\right].
\end{eqnarray}
Note that $V(-\delta_1,\ldots,-\delta_N)$ is an invertible
matrix because, by assumption, for all $i=1,\ldots,N, j=1,\ldots,N,$
if $i\not=j$ then
$\delta_i\not=\delta_j$ and therefore we have
$\det\left(V(-\delta_1,\ldots,-\delta_N)\right)=
\Pi_{1\leq i<j \leq N}(\delta_i-\delta_j)\not=0$
(cf. e.g. \cite{LancTism}, p.35).

The first order equations can now be rewritten as
\begin{eqnarray}
\label{focVI}
\left[\begin{array}{c}
\tilde{a}(-\delta_1)^2\\
\vdots\\
\tilde{a}(-\delta_{N})^2\end{array}\right] & = &
\mbox{diag}(e(\delta_1),\ldots,e(\delta_N))V(\delta_1,\ldots,
\delta_N)
V(-\delta_1,\ldots,-\delta_N)^{-1}
\left[\begin{array}{c}
\tilde{a}(-\delta_1)\\
\vdots\\
\tilde{a}(-\delta_N)
\end{array}\right]\nonumber
\end{eqnarray}
\begin{equation}
[\tilde{a}(-\delta_1),\ldots,\tilde{a}(-\delta_N)] \not= 0
\end{equation}
where $\mbox{diag}(e(\delta_1),
\ldots,e(\delta_N))$ denotes the diagonal
matrix with $e(\delta_i)$ in the $(i,i)-$entry, $i=1,\ldots,N.$\\
This means that these first order equations
can be written as
\begin{eqnarray}
\label{focVII}
\left[\begin{array}{c}
      x_1^2\\
      x_2^2\\
      \vdots\\
      x_N^2
      \end{array}\right]&=&M
\left[\begin{array}{c}
      x_1\\
      x_2\\
      \vdots\\
      x_N
      \end{array}\right],~x \not=0
\end{eqnarray}
where $x_i=\tilde{a}(-\delta_i),
i=1,\ldots,N,~x=(x_1,\ldots,x_N)'$
and
\begin{equation}
\label{M}
M=\mbox{diag}(e(\delta_1),\ldots,e(\delta_N))
V(\delta_1,\ldots,
\delta_N)
V(-\delta_1,\ldots,-\delta_N)^{-1}.
\end{equation}
In the next section the solution of equations of the form found
here will be treated in general.

\section{Diagonal-quadratic systems of equations}
\label{quaddiag}

In this section  we will present results about an arbitrary
system of polynomial equations of the form
\begin{eqnarray}
\label{diagquadeq}
\left[\begin{array}{c}
      x_1^2\\
      x_2^2\\
      \vdots\\
      x_N^2
      \end{array}\right]&=&M
\left[\begin{array}{c}
      x_1\\
      x_2\\
      \vdots\\
      x_N
      \end{array}\right]+\mu,
\end{eqnarray}
where $\mu \in {\bf C}^N$ is a constant $N-$vector.
This will be called a {\em diagonal-quadratic} system of equations.

{\em Remark.} A quadratic equation in $x$ can be written
as $x^TAx+cx+d$ for some symmetric matrix $A,$ a row vector $c$ and
 a scalar $d.$ If $A=e_i e_i^T,$ for some $i \in \{1,\ldots,N\},$
then the equation is one of the form described above. If there are
$N$ quadratic equations and the corresponding
$A-$matrices are all diagonal, and these diagonal matrices form
a basis of the linear vector space of all diagonal $N \times N$
matrices then such a system can (obviously) be rewritten in the
form above. That is the motivation for the terminology
`diagonal-quadratic'.

In this paper use will be made of Gr\"{o}bner basis theory and
constructive algebra. For an exposition of this theory one can
refer to e.g. \cite{CLOS}. In Gr\"{o}bner basis theory an
important
role is played by the so-called {\em monomial orderings}. Let
$\alpha=(\alpha_1,\ldots,\alpha_N)$ denote an arbitrary vector
of nonnegative integers, which will be called a multi-index in the
sequel, then $x^{\alpha}$ will denote the monomial
$x^{\alpha}:=x_1^{\alpha_1}x_2^{\alpha_2}\ldots x_N^{\alpha_N}.$
The multi-index $\alpha$ is called the {\em multi-degree}  of
the monomial $x^{\alpha}.$  The corresponding total degree is
defined as $|\alpha|:=\alpha_1+\alpha_2+\ldots+\alpha_N.$ For a
general definition of monomial ordering we refer
to \cite{CLOS}, p.54, Definition 1.

A partial ordering of monomials is defined by $x^{\alpha} \succ
x^{\beta}$ if $|\alpha|>|\beta|$. Such an ordering is called a total
degree ordering. For our purposes any complete ordering which is a
refinement of the total degree ordering would do.
For definiteness we choose to work with the graded lexicographic
ordering, which refines the total degree ordering as follows: if
$|\alpha|=|\beta|$ then $x^{\alpha} \succ x^{\beta}$ if $\alpha_i>
\beta_i$ for the smallest integer $i \in \{1,\ldots,N\}$ for which
$\alpha_i \not= \beta_i$.

The total degree of a polynomial
is defined as follows. Each polynomial is a unique linear combination of
monomials
with nonzero coefficients. The maximal  total degree of
these monomials is called the {\em total degree of the
polynomial}.
If we denote the $i-$th row of the
matrix
$M$ by $m_i$ and the $i-$th entry of the vector $\mu$ by $\mu_i$
for $i \in \{1,\ldots,N\},$ then the equations
can be rewritten as
\[x_i^2-m_i x -\mu_i=0, i=1,\ldots,N.\]
Let $g_i(x_1,\ldots,x_N):=x_i^2-m_i x -\mu_i, i=1,\ldots,N,$
 then we are looking for the zeros of the ideal $I$ spanned by
$G:=\{g_1,g_2,\ldots,g_N\}.$

Let $<g_1,\ldots,g_N>$ denote the ideal generated by the set of
polynomials $g_1,\ldots,g_N$. For a polynomial
$f$, let $LT(f)$ denote the leading term of $f$, and for an ideal $I$
of polynomials, let $LT(I)$ denote the set of all leading terms of the
polynomials in $I$.

\begin{definition}
For a fixed monomial ordering, a finite subset
$\Gamma=\{\gamma_1,\ldots,\gamma_\nu\}$ of
an  ideal $I$ is a Gr\"{o}bner basis if
\[ <LT(\gamma_1),\ldots,LT(\gamma_\nu)> = <LT(I)>. \]
\end{definition}

\begin{theorem}
\label{quaddiag:groebner}
The set $G$ is a Gr\"{o}bner basis with respect to total degree ordering.
\end{theorem}
{\bf Proof.}
With respect to any ordering which is a refinement of partial ordering
by total degree,
the leading terms of $G$ are monomials of the form $x_i^2$. These are
clearly pairwise coprime. But it is known that this implies that $G$
is a Gr\"{o}bner basis~\cite[p.333, Ex.15.20]{Eisenbud}.
\phantom{xxxxx} \hfill $\square$\\
An  alternative but longer proof is available in \cite{HanMacChoTR314}.

This result is very important because to apply the results of
Gr\"{o}bner basis theory one needs a Gr\"{o}bner basis. Usually
one needs to apply an algorithm like Buchberger's algorithm to
bring a set of polynomials that generates the ideal in which one
is interested in Gr\"{o}bner basis form. In fact in many cases
this is the most difficult part of the calculations. In the case
at hand however the set of polynomials of which we want to find
the zeros itself forms a Gr\"{o}bner basis.

But that is not all. We can say more. We know that
$G=\{g_1,\ldots,g_N\}$ forms a Gr\"{o}bner basis and that the leading
monomial of $g_i$ is $x_i^2$ for each $i=1,\ldots,N.$
Let ${\bf C}[x_1,\ldots,x_N]$ denote the
ring of polynomials with complex coefficients.
Let $R$ denote the set of multi-indices $R:=\{0,1\}^N.$
In other words, $R$ is the set of all multi-indices
$\alpha=(\alpha_1,\ldots,\alpha_N)$ with the property that
for each $i=1,\ldots,N$ one has either
$\alpha_i=0$ or $\alpha_i=1.$ Let $Q$ denote the set of all
multi-indices outside $R.$ For each polynomial $p=p(x)$ there
exists
a unique additive decomposition $p=p^R+p^Q,$ where the polynomial
$p^R$ is a linear combination of monomials with multi-degree in
$R$ and $p^Q$ is a linear combination  of monomials with
multi-degree in $Q.$
\begin{lemma}
\label{Sbasis}
Let $I$ denote the ideal generated by $G.$
\begin{itemize}
\item[(i)] The set $V=V(I)$ of zeros in ${\bf C}^N$
of the system
of polynomial equations $g_i(x)=0,~i=1,\ldots,N,$
is finite.
\item[(ii)] The ${\bf C}-$vector space $S=Span(x^{\alpha}:x^{\alpha} \not\in
<LT(I)>)$ is finite-dimensional.
\item[(iii)] The ${\bf C}-$vector space ${\bf C}[x_1,\ldots,x_N]/I$
is finite-dimensional.
\item[(iv)]
The set of monomials $\{x^{\alpha}:\alpha \in R\}$ forms a basis
for the vector space $S$.
\item[(v)]
The dimension of the vector space $S$ is $2^N.$
\item[(vi)]
The dimension of the vector space ${\bf C}[x_1,\ldots,x_N]/I$
is $2^N.$
\end{itemize}
\end{lemma}

{\bf Proof.}
ad(i)--(iii). (i)---(iii) follow immediately from \cite[Chapter 5,
Theorem 6]{CLOS}.\\
ad (iv). Because $G$ is a Gr\"{o}bner basis the ideal
$<\mbox{LT}(I)>$ is equal to the ideal generated by the leading
terms of the elements of $G,$ i.e. the ideal $<x_1^2,\ldots,x_N^2>.$ The
monomials in this ideal are precisely  those which have multi-degree in the
set $Q.$ Therefore the monomials
in $S$ are the all the monomials with multi-degree in $R.$
\\
ad (v). From (iv) it follows that the dimension of $S$ is equal
to the cardinality of $R, $ which is $card(R)=2^N.$\\
ad (vi). According to Proposition 4 of Chapter 5 of \cite{CLOS}
the vector space
 ${\bf C}[x_1,\ldots,x_N]/I$ is isomorphic
to $S$ and therefore has the same dimension as $S.$
\phantom{xxxxx} \hfill $\square$

From \cite{CLOS}, Chapter 5, Section 3, Proposition 1 it follows that
every polynomial in ${\bf C}[x_1,\ldots,x_N]$ can be written in a
unique way as the sum of an element of $S$ and an element of $I.$ In
other words, each equivalence class $f+I,$ where $f$ is an arbitrary
polynomial in ${\bf C}[x_1,\ldots,x_N],$ has a unique representative
in $S.$ Let this representative be denoted by $\pi(f) \in S.$ Given
$f,$ the polynomial $\pi(f)$ can be obtained by a general method from
Gr\"{o}bner basis theory, namely the so-called division algorithm with
respect to the Gr\"{o}bner basis $G$ as described in e.g. \cite{CLOS}.
However, for diagonal quadratic equations, the division algorithm
simplifies considerably and one can describe in direct terms how one
can obtain $\pi(f)$ from $f.$ The `reduction procedure' can be
described as follows.  Using the same notation as above, one can write
$f=f^Q+f^R,$ where $f^R \in S$ and the monomials of $f^Q$ all have
multi-degree in $Q.$ This additive decomposition is obviously unique.
If $f^Q=0$ then $f=f^R \in S$ in which case $\pi(f)=f$ and we are
done. If $f^Q\not=0$ then consider any monomial of $f^Q$ with total
degree equal to the total degree of $f^Q.$ By construction each such
monomial is divisible by at least one of the monomials
$x_1^2,x_2^2,\ldots,x_N^2.$ If it is divisible by $x_i^2$ then
replacing it by the polynomial that is obtained by multiplying the
monomial by $\frac{h_i(x)}{x_i^2}$ the result is a polynomial
$\tilde{f}$ that is in the equivalence class $f+I$ and which has the
following property. Either the total degree of $\tilde{f}^Q$ is
smaller than the total degree of $f^Q,$ or otherwise the total degree
of $\tilde{f}^Q$ is equal to the total degree of $f^Q$ but the number
of monomials in $\tilde{f}^Q$ with total degree equal to the total
degree of $f^Q$ is one less than the number of monomials in $f^Q$ with
total degree equal to the total degree of $f^Q.$ Such a replacement of
$f$ by $\tilde{f}$ will be called a `reduction step'. It follows that
after a finite number of such reduction steps one arrives at a
polynomial in the equivalence class $f+I$ with the property that it
lies in $S.$ This is then the unique polynomial $\pi(f)$ that was
sought for.

The importance of this reduction procedure in our application
 will become clear in
 the
examples section.

\section{Commutative matrix solutions of polynomial equations}
In this section a method to obtain the solutions of a system of
polynomial equations in several variables will be outlined. A
method of this kind was originally developed by
\cite{Stetter,MollerStetter}. A similar approach, but differing in
some details, was developed by the authors of the present paper,
is available in \cite{HanMacChoTR314}, and is the approach which
will be summarized here. All proofs are omitted from this section
since they are available in the works cited above.

We will consider the situation in which the system of polynomial
equations will have a finite number of solutions over the field of
complex numbers {\bf C}. In the modern constructive algebra
approach to the problem of finding the roots of a system of
polynomial equations the theory of Gr\"{o}bner bases plays an
important role. For this theory we refer, as before, to
\cite{CLOS}. A fundamental theorem of the theory of Gr\"{o}bner
bases is that for any polynomial ideal given by a finite number of
polynomials which generate it, a Gr\"{o}bner basis can be
calculated with respect to any admissible
 monomial ordering
(like the lexicographical ordering or the total degree ordering)
in a {\em finite} number of steps. It can for example be obtained
by Buchberger's algorithm. However the number of steps required by
such an algorithm can be huge.
In the literature
it is suggested that in order to obtain the roots of a system of
polynomial equations, one can  construct a Gr\"{o}bner
basis with respect to a lexicographical ordering
\cite[p.233]{CLOS},\cite[pp. 459-462]{GeddesCzaporLabahn}.
Also in the paper \cite{HanMacAut96} this approach was followed
to show that under two hypotheses described in that paper, the
 $H_2$ model order reduction
problem can be solved in a finite number of steps.
However only examples of
reduction of third order models were presented in
 that paper. The bottle-neck
in the calculations was the construction of a Gr\"{o}bner basis.
In the previous section it was shown that for the problem of
reduction of the model order by one with respect to the $H_2$
norm, in case of an original model with distinct poles, the first
order equations found already are in the form of a total degree
Gr\"{o}bner basis, so
 {\em no Gr\"{o}bner basis construction
at all} is required in the application at hand.

The idea is first to construct a {\em commutative matrix solution} for
a system of polynomial equations which is in Gr\"{o}bner basis
form.
\begin{definition}
Let $N$ be a positive integer.
Let $f \in {\bf C}[x_1,\ldots,x_N]$ be a polynomial
in the variables $x_1,\ldots,x_N.$
Let $M$ be a positive integer and consider
an $N-$tuple $(A_1,A_2,\ldots,A_N)$ of square $M \times M$
matrices that commute with each other, i.e. $A_iA_j=A_jA_i$
for each pair $(i,j), i=1,\ldots,N,~j=1,\ldots,N.$
Then $(A_1,A_2,\ldots,A_N)$ will be called a
commutative matrix solution of the  polynomial equation
$f=0$ if $f(A_1,\ldots,A_N)=0_M,$ where the symbol $0_M$ denotes
the $M \times M$ zero matrix.
\end{definition}
In the following, an $M\times M$ zero matrix will often be denoted
by the symbol $0,$ as is usual, instead of the symbol $0_M.$ The
size of the zero matrix should then be clear from the context.
An $N-$tuple of $M \times M$ matrices $(A_1,\ldots,A_N)$
will be  called a commutative matrix solution of a system of
polynomial
equations in $N$ unknowns $x_1,\ldots,x_N,$ if it is a commutative
matrix solution for each of the polynomials in the system.

From a commutative matrix solution a scalar solution
can be obtained by considering any common eigenvector
of the matrices. The corresponding eigenvalues form
an $N-$tuple which is in fact a scalar solution of
the system of polynomial equations involved.
 The commutative matrix solution that will be constructed
here for the case of ideals with zero dimensional variety, has the
 property that ALL (scalar) solutions can be obtained
in this way.

It will first be explained how such a commutative matrix solution
can be constructed. Then it will be shown how the (scalar)
 solutions
of the system of polynomial equations can be derived from the
 matrix solution by eigenvalue-eigenvector calculations.
If ${\cal F}$ is a field containing all the coefficients
of the polynomial equations then all the entries of the matrix
 solution that will be constructed
 will be contained in ${\cal F}$; in other words,
only additions, subtractions, multiplications and divisions are
required to obtain a matrix solution.

We start with two results which hold for an arbitrary
polynomial ideal. For these results to hold, the ideal
  does {\em not} have to have the property
that the number of common zeros of the polynomials in the
 ideal is finite.
The two results consist of  a number of observations concerning
the operation `multiplication by $x_i$ modulo the ideal',
for $i \in
\{1,\ldots,N\}.$
Composition of a pair of mappings $X,~Y$ will be denoted
 (as usual) by $X \circ Y.$

\begin{theorem}
Let $N$ be a positive integer.
Let $I \subset {\bf C}[x_1,\ldots,x_N]$ be an ideal and let
${\cal V}:= {\bf C}[x_1,\ldots,x_N]/I$ denote the corresponding
residue class ring.
Let $i \in \{1,\ldots,N\}$ be fixed.
 Let $f_1,f_2 \in {\bf C}[x_1,\ldots,x_N].$
If $f_1$ and $f_2$ are equal modulo $I,$
then  $x_i f_1$ and $x_i f_2$ are equal modulo $I.$
The mapping $X_i:{\cal V} \rightarrow {\cal V}, f+I \mapsto x_i
f+I,$
is a linear endomorphism.
For $i, j \in \{1,\ldots,N\}$ arbitrary, $X_i \circ X_j=X_j \circ X_i$
i.e. the linear mappings $X_i$ and $X_j$ commute. The mapping $X_i
\circ X_j$ is the mapping
given by $f+I \mapsto x_i x_j f+I.$
\end{theorem}

For any pair of linear endomorphisms $X,~Y$ let us interpret $XY$
as the composition $X \circ Y,$ let us interpret $X^0$ as the
identity  and for each positive integer $k,$ let us interpret the
power $X^k$ as the $k-$fold composition $X \circ X \circ \ldots
\circ X.$ Using this interpretation for any $N-$tuple of {\em
commutative} linear endomorphisms $X_1,\ldots,X_N$ and any
polynomial $f \in {\bf C}[x_1,\ldots,x_N], $ the expression
$f(X_1,\ldots,X_N)$ denotes a well-defined linear endomorphism.

\begin{theorem}
Let $N, I, {\cal V}$ and $X_i,~i=1,\ldots,N$ be as given in
 the previous theorem.
For any polynomial $f \in {\bf C}[x_1,\ldots,x_N]$ the linear
mapping $f(X_1,X_2,\ldots,X_N): {\cal V} \mapsto {\cal V}$ is
well-defined.

The following two  statements are equivalent,
\begin{itemize}
\item[(i)] $f \in I,$
\item[(ii)] $f(X_1,\ldots,X_N)$ is equal to the zero mapping
${\cal V} \rightarrow {\cal V}, f+I \mapsto 0+I.$
\end{itemize}
\end{theorem}

Now we will specialize to systems of polynomial equations with
finitely many common solutions.
We will make extensive use of the results from section 3 of
Chapter 5 of \cite{CLOS}, pp. 228-235, especially Propositions 1
and 4 and Theorem 6 of that section.

Let $g_1(x_1,\ldots,x_N)=0,\ldots,g_{N'}(x_1,\ldots,x_N)=0$ denote
a system of $N'$ polynomial equations with complex coefficients
 in the $N$
 variables $x_1,\ldots,x_N.$ The complex vector
 $(\xi_1,\ldots,\xi_N) \in {\bf C}^{N}$ is a root of the system
 of polynomial
equations
if for each $j=1,\ldots,N',$
\[g_j(\xi_1,\ldots,\xi_N)=0.\]
Let $I=<g_1,\ldots,g_{N'}> \subset {\bf C}[x_1,\ldots,x_N]$ denote
 the ideal generated  by the polynomials
$g_1(x_1,\ldots,x_N),$ $\ldots, g_{N'}(x_1,\ldots,x_N).$

Suppose that $G=\{g_1,\ldots,g_{N'}\}$ is in fact a Gr\"{o}bner
 basis for $I$, with respect to some fixed monomial ordering.
Similarly to what was noted in the previous section for the
special case of diagonal-quadratic systems of polynomial
equations, the following can be said for this more general case.
Each polynomial  $f \in {\bf C}[x_1,\ldots,x_N]$ is congruent
modulo $I$  to a polynomial $r$ with leading term
that cannot be reduced by any of the leading terms of the
polynomials in the Gr\"{o}bner basis; for each $f$ the associated
polynomial $r$ is unique \cite[Chapter 5, Section 3, Proposition
1]{CLOS} and will be denoted by $\overline{f}^G.$ The set $V$ of
all polynomials $r$ obtained in this way forms
 a finite dimensional vector space if and only if
the number of roots of the system of polynomial equations is
finite. If this set is indeed a finite dimensional vector space,
then it has a basis consisting of monomials, namely all monomials
that cannot be reduced by any of the leading terms of the
polynomials in the Gr\"{o}bner basis. This result is due to
Macaulay \cite[Theorem 15.3, p.325]{Eisenbud}. Given the monomial
ordering it is a straightforward task to list these monomials
(\cite{CLOS}). Let this basis be denoted by $B.$
The mapping $V \rightarrow {\cal V}, r \mapsto r+I,$
is a linear bijection of vector spaces. In case ${\cal V}$ is
 finite dimensional, let ${\cal B}$ denote
the basis of ${\cal V}$ obtained as the image of $B$ under this
mapping. Let $D$ denote the dimension of ${\cal V}.$

For each $i \in \{1,\ldots,N\}$ let
$A_{X_i}$ denote the $D \times D-$matrix of the endomorphism
$X_i$ with respect to the basis ${\cal B}.$

Using this set-up the following fundamental result can be
obtained.
\begin{theorem} \label{thm:fAXizero}
Let a monomial ordering be fixed and let $G$ be a Gr\"{o}bner
 basis of the ideal $I.$
Let the associated linear space ${\cal V}$ be finite dimensional
 with dimension $D.$
Let $f \in {\bf C}[x_1,\ldots,x_N]$ be given. Let the mappings
 $X_i,~i=1,\ldots,N$ and $f(X_1,X_2,\ldots,X_N)$ be as given in
 the previous theorems.

 The matrix of the linear mapping
$f(X_1,X_2,\ldots,X_N):{\cal V} \rightarrow {\cal V}$ with respect to
the basis of monomials ${\cal B}$ of ${\cal V}$ is equal
to $f(A_{X_1},A_{X_2},\ldots,A_{X_N}).$

The following two statements are equivalent,
\begin{itemize}
\item[(i)] $f \in I,$
\item[(ii)] $f(A_{X_1},A_{X_2},\ldots,A_{X_N})=0,$
i.e. this matrix is the $D \times D$ zero matrix.
\end{itemize}
\end{theorem}

This theorem tells us that the $N-$tuple of matrices
$(A_{X_1},\ldots,A_{X_N})$ is in fact a commutative matrix
 solution of any
system of polynomial equations that generates $I.$

The entries  of the $k-$th column of the matrix $A_{X_i}$ are
obtained as follows. Let the $k-$th element of the basis $B$ of
monomials be denoted by $b_k.$ The monomial $x_i b_k$ is either
itself in the basis $B$ or otherwise
 $\overline{x_i b_k}^G\not= x_i b_k.$ In both
cases  $\overline{x_i b_k}^G$ can be written as a unique linear
combination of the elements of $B.$ The coefficients of the linear
combination are the entries of the $k-$th column of  the matrix
$A_{X_i}.$ If $x_i b_k$ is itself in the basis $B$ then the $k-$th
column of the matrix $A_{X_i}$ is a standard basis vector.

In the case $N=1$ then there exists a unique monic polynomial $g$ such
that $I$ is generated by $g.$ In that case the matrix $A_{X_1}$ is a
{\em companion matrix} of $g$ (cf. e.g. \cite[p. 68]{LancTism}).

Now suppose that the vector $v$ is a {\em common} eigenvector
of the matrices $A_{X_1},\ldots,A_{X_N}$ with corresponding
eigenvalues $\xi_1,\xi_2,\ldots,\xi_N,$ respectively, i.e.
for each  $i \in \{1,\ldots,N\}$ the equality
$A_{X_i} v=\xi_i v$ holds and $v \not= 0.$
Then for each $f \in I$ one has
\[0=f(A_{X_1},\ldots,A_{X_N})v=f(\xi_1,\ldots,\xi_N)v\]
and therefore $f(\xi_1,\ldots,\xi_N)=0.$ It
follows that
$(x_1,\ldots,x_N)=(\xi_1,\ldots,\xi_N)$ is a
root of any system of polynomial equations that generates the ideal
$I.$

The following fundamental result states that in fact
{\em each} of the finite number of roots is obtained in
this way.
\begin{theorem}
Let $N$ be a positive integer and let $I$ be an ideal in
the ring ${\bf C}[x_1,\ldots,x_N]$ such that the corresponding
set $Z \subset {\bf C}^N$ of common zeros of all the
polynomials in $I$ is finite. Let $X_i,~i=1,\ldots,N$ be
as defined above.
Then for each vector $\xi=(\xi_1,\ldots,\xi_N)'
\in Z$ there exists a polynomial  $w \in {\bf C}[x_1,\ldots,x_N],$
$w \not\in I,$ with the property that for each $i=1,\ldots,N,$
the following equality holds:
\[X_i(w+I)=\xi_i w+I,\]
i.e. $w$ is a common eigenvector of the mappings
$X_1,X_2,\ldots,X_N,$ with corresponding eigenvalues
$\xi_1,\ldots,\xi_N,$ respectively.
\end{theorem}

From this theorem we have the following important corollary.
\begin{corollary}
\label{eigvalsol}
Let $N,$ $I$ and $Z$ be as given in the previous theorem.
Let $X_i,~i=1,\ldots,N$ be as defined above.
Let a monomial ordering be given and let $G$ be a
Gr\"{o}bner basis of $I$ with respect to this monomial
ordering. Let $B$
denote the basis of all monomials in ${\bf C}[x_1,\ldots,x_N]$ that
are not included in the ideal $<LT(G)>$ generated by the leading
terms of the elements of $G$ and let ${\cal B}$ denote the
corresponding
basis of the residue class ring ${\bf C}[x_1,\ldots,x_N]/I,$
as before. Let
$A_{X_1},\ldots,A_{X_N}$ denote the matrices of the linear
endomorphisms $X_1,\ldots,X_N, $ respectively, with respect
to the basis
${\cal B}.$
Then the following two statements are equivalent.
\begin{itemize}
\item[(i)] $\xi=(\xi_1,\ldots,\xi_N)' \in Z.$
\item[(ii)] There exists a common eigenvector
$v \in {\bf C}^N \setminus \{0\}$ of the (commutative)
matrices $A_{X_1},\ldots,A_{X_N}$ with corresponding eigenvalues
$\xi_1,\ldots,\xi_N$ respectively, i.e. there exists a
nonzero vector $v$ with the property
\[A_{X_i}v=\xi_i v,~i=1,\ldots,N.\]
\end{itemize}
\end{corollary}

Various alternatives arise as to how to exploit the theory
presented here to solve a system of polynomial equations, starting
with a Gr\"{o}bner basis. The commutative matrix solution
presented can be calculated in symbolic form if the original
system of equations is in symbolic form
and it can be calculated exactly in numerical form if the
coefficients of the original system of polynomials are given
numerically. From the commutative matrix solution the roots of the
system of polynomial equations can be obtained either by exact
algebraic methods or by numerical methods that involve round-off
errors. The exact algebraic approach will not be worked out here.

A (nonexact) numerical approach can be based on numerical
calculation of the eigenvalues and eigenvectors of the matrices
involved. In the examples section this approach will be applied to
the $H_2-$model order reduction problem.

The possibility of using a mixture of exact and symbolic
calculations with numerical calculations is  very promising for
obtaining practically useful results. The matrices involved will
tend to become huge (in terms of numbers of entries) if the number
of variables involved grows; however eigenvalue calculation can be
done numerically for quite big matrices. In section
\ref{sec:examples} matrices with several hundreds of rows and
columns are used. One can expect that usage of more
 refined numerical
techniques will make it possible to push the limits
quite a bit further.

Let $f \in {\bf C}[x_1,\ldots,x_N]$ and let $F$ be the
corresponding linear endomorphism of ${\bf C}[x_1,\ldots,x_N]/I$
defined by $g+I \mapsto f.g+I.$ If the number of common zeros of
the polynomials in $I$ is finite, and we have a basis ${\cal B}$
of
  ${\bf C}[x_1,\ldots,x_N]/I$
as before, then we can represent $F$ with respect to this basis by
a matrix $A_F.$ It is now straightforward to see that
$A_F=f(A_{X_1},A_{X_2}, \ldots, A_{X_N}).$ More generally if
$f=\frac{f_n}{f_d},~f_n,f_d \in {\bf C}[x_1,\ldots,x_N]$ and
$f_d(\xi)\not=0$ for each common zero $\xi$ of the polynomials in
$I,$ then $F$ and $A_F$ are again well-defined and
$A_F=f_n(A_{X_1},\ldots,A_{X_N}).
\left(f_d(A_{X_1},\ldots,A_{X_N})\right)^{-1}.$ The eigenvalues of
this matrix $A_F$ are $\{f(\xi)|\xi\in Z\}$. For example in
optimization problems in which the criterion function $f,$ say, is
a rational function this can be used to obtain the matrix $A_F$
which has as its eigenvalues the critical values of $f$. (The
values that a function takes on its set of critical points are
called the critical values.) The matrix $A_F$ could be called a
{\em critical value matrix} and its characteristic polynomial a
{\em critical value polynomial}.
 This is related to Theorem 9 and the subsequent
Remark 10 in \cite{HanMacAut96} concerning the existence
and usage of a univariate polynomial which has the critical values
of the criterion function as its zeros.

\section{Model order reduction by one in $H_2$}
\label{red1sol}

Recall the formulation of the $H_2$ model reduction problem from
Section \ref{foc}. In order to facilitate the statement of the
following theorem let us define the set $\Xi$ as follows.  Let
$\frac{e}{d} \in \Sigma S_{N}$ have $N$ distinct poles
$\delta_1,\ldots,\delta_N \in {\bf C}.$ Let the matrix $M$ be as
given in equation (\ref{M}) and let $\Xi$ denote the set of
solutions in ${\bf C}^N \setminus \{0\}$ of equation
(\ref{focVII}). The diagonal quadratic system of equations
(\ref{focVII}) is shown to form  a total degree Gr\"{o}bner basis
in Theorem \ref{quaddiag:groebner}. In Lemma \ref{Sbasis} a basis
of $2^N$ monomials of the corresponding vector space $S$ is
presented. This basis consists of the monomials outside the ideal
generated by the leading terms of all polynomials in the ideal
corresponding to the diagonal quadratic equations. Let this basis
be denoted by $B.$  Then Corollary \ref{eigvalsol} can be applied
to (\ref{focVII}) using the basis of monomials $B.$ The
implication is that in this case the set $\Xi$ just defined is
equal to the set $Z$ of that Corollary, except that the zero
vector is removed:
\[\Xi=Z \setminus \{0\} \]
It follows that $\Xi$ contains at most $2^N-1$ elements,
each of which is a vector of $N$ entries that can be found as the
eigenvalues
corresponding to  any common eigenvector
of the matrices $A_{X_1},\ldots,A_{X_N}$ from Corollary \ref{eigvalsol}.
We
therefore have the following theorem

\begin{theorem}
\label{th:H2critical}
Let $\frac{e}{d} \in \Sigma S_{N}$ have $N$ distinct poles
$\delta_1,\ldots,\delta_N \in {\bf C}.$
\begin{itemize}
\item[(i)] The number of critical points of the criterion function
$f: \Sigma S_{N-1} \rightarrow [0,\infty), \frac{b}{a} \mapsto
\left\|\frac{e}{d}-\frac{b}{a}\right\|_2^2$ is finite and not
greater than $2^N-1.$ \item[(ii)] If the rational function
$\frac{b}{a} \in \Sigma S_{N-1}$ is a critical point of $f$ then
there exists a number $q_0$ and a vector $\xi \in \Xi \subset {\bf
C}^N \setminus \{0\}$ such that
 $q_0 a(-\delta_i)=\xi_i,~i=1,\ldots,N.$
For given $q_0$ and $\xi$ the polynomial $a$ is uniquely
 determined by this linear system of equations and $b$ is
 uniquely determined by equation (\ref{bequation}).
\end{itemize}
\end{theorem}

Of course the solutions that will be found for the first order
equations will in general not all correspond to rational functions
$\frac{b}{a} \in \Sigma S_{N-1}$: it is certainly possible that
some
will not correspond to real systems; some may correspond to real
but
unstable systems.

An algorithm to obtain all the critical points of the criterion
function of $H_2$ model reduction by one could now be constructed
as follows.
\begin{enumerate}
\item Construct the matrix $M.$
\item Construct the matrices $A_{X_1},\ldots,A_{X_N}.$
\item Calculate the eigenvalues of these matrices that correspond
to a common eigenvector of all these matrices. The result is
a vector $\xi \in {\bf C}^N.$ All nonzero vectors $\xi$
obtained in this way form the (finite) set $\Xi.$
\item For each element of $\Xi$ solve equation
(\ref{focVII}) for $a$ and $q_0$, and
select those $a$ that are real and Hurwitz.
\item For those $a$ selected in the previous step, solve equation
  (\ref{bequation}) for $b$.
\end{enumerate}

Note that steps (1) and (2) can be done purely symbolically. Apart
from considerations of memory storage and perhaps calculation
time, it is not necessary to specify the original system; one can
present it symbolically by its poles and the (non-zero) values of
the numerator polynomial in these poles.

If the original system is specified numerically then step (3) can
be worked out by either constructive algebra algorithms (using e.g.
methods of isolation of zeros of polynomials) or by numerical
algorithms that admit round-off errors.  In section
\ref{sec:examples} we present some results obtained by
calculations of the latter type.

Step (4) requires that we go through the solutions in $\Xi$ to find
out those that are admissible and a solution is admissible if $a$ is
both real and Hurwitz. This can be done by first eliminating all the
complex $a$'s and then checking whether the real $a$'s are Hurwitz.

Note that the pairs $a,b$ found in Steps (4) and (5), respectively, are
coprime as a consequence of equation (\ref{bequation}), and that $b$
is real, and hence that $\frac{b}{a} \in \Sigma S_{N-1}.$

The global approximant is found by selecting from the finite set of
critical points the point that minimizes the criterion function $f$
defined in Theorem \ref{th:H2critical}. This follows from the fact
that this criterion function $f$ is differentiable everywhere and has
a global minimum (cf. \cite{baratchartolivi88} and the references
therein).  The global approximant can now be found by choosing the
admissible solution that gives the least $H_2$ criterion function.

This process can be simplified, since one is interested in
locating only the global approximant. We shall show that it is
possible to construct a matrix, the eigenvalues of which include the
values of the criterion function $f$ at the critical points. One can
therefore search among these values, starting with the smallest
positive real value, until one finds one which corresponds to an
admissible approximant. This will be the optimal approximant. As will
shortly be shown, the
attraction of this  approach is that many elements of $\Xi$, namely
those which yield complex value of $f$ and those which correspond to
non-Hurwitz $a$ polynomials, will never be visited by this procedure.

For any rational function $t$ for which the Lebesgue integral
$\frac{1}{2 \pi} \int_{-\infty}^{\infty} |t(i\omega)|^2 d\omega$
is finite let us define the $L_2-$norm $\|t\|_2$ by
\[\|t\|_2^2:=\frac{1}{2 \pi}
\int_{-\infty}^{\infty} |t(i\omega)|^2 d\omega.\]
Note that for any rational function $t$ in $H_2$ this definition
coincides with the definition of $\|t\|_2$ given before.
We have the following theorem.
\begin{theorem}
\label{realsol} Let $\frac{e}{d} \in \Sigma S_{N}$ have $N$
distinct poles $\delta_1,\ldots,\delta_N \in {\bf C}.$

Let $a(s), b(s), q_0$ be a real solution of  the polynomial
equations (\ref{bequation}),(\ref{focIII}), then
\begin{enumerate}
\item
\begin{eqnarray}
\left\| \frac{e(s)}{d(s)} - \frac{b(s)}{a(s)} \right\|_2^2 & = &
\sum_{i = 1}^N \frac{x_i^3}{e(\delta_i)d'(\delta_i)d(-\delta_i)}
\label{thirddegree}
\end{eqnarray}
where $x_i=\tilde{a}(-\delta_i) = q_0 a(-\delta_i),~i=1,\ldots,N,$
(as before) and $d'(s)$ denotes the derivative of
$d(s)$ with respect to $s$.
\item If $a(s)$ is Hurwitz then the $L_2$-norm computed above
coincides
  with the $H_2$-norm.
\item If $a(s)$ is not Hurwitz then the $L_2$-norm squared
computed above is
  strictly greater than the global minimum of the criterion
function $f$ as defined in Theorem \ref{th:H2critical}.
\end{enumerate}
\end{theorem}

{\bf Proof.}
Let us first prove part 1 of the theorem. Due to the first order
condition (\ref{focII}), combined with
the equality $q(s)=q_0$ and combined with the assumption that
$e,d,a,b$ are real polynomials, and combined with the fact that
$d$ and $a$ are monic polynomials and therefore unequal to the
zero polynomial, one has
\[\left\| \frac{e(s)}{d(s)} - \frac{b(s)}{a(s)} \right\|_2^2 =
\left\| \frac{a(-s)^2 q_0}{a(s)d(s)} \right\|_2^2=\]
\[\frac{1}{2 \pi}
\int_{-\infty}^{\infty} \frac{a(-i\omega)^2 a(i\omega)^2
q_0^2}{d(i\omega) a(i\omega) d(-i\omega) a(-i\omega)} d\omega=
\frac{1}{2 \pi}\int_{-\infty}^{\infty} \frac{a(-i\omega) a(i\omega)
q_0^2}{d(i\omega) d(-i\omega)} d\omega.\]
The residue theorem of complex analysis can now be applied.
We use the fact that $\lim_{|s|\rightarrow \infty} s^2
 \left(\frac{a(-s) a(s)}{d(s) d(-s)}\right)=1$ to argue that
the integral over the imaginary axis is equal to the integral
over a sufficiently large semi-circle together with a sufficiently
large segment of the imaginary axis. This is a standard argument
 in complex analysis that we will not repeat here (see e.g.
\cite{NevanPaat}). The residue theorem now
tells us that the
integral is equal to
\[ q_0^2 \sum_{i=1}^{N} \mbox{Res}_{s=\delta_i} \left(\frac{a(-s)
a(s)}{d(s) d(-s)}\right)= \]
\[ q_0^2 \sum_{i=1}^{N} \lim_{s \rightarrow \delta_i}
\left(\frac{(s-\delta_i) a(-s) a(s)}{d(s) d(-s)}\right)= \]
\[ q_0^2 \sum_{i=1}^{N} \left(\frac{a(-\delta_i)
a(\delta_i)}{d'(\delta_i) d(-\delta_i)}\right)= \]
\[ \sum_{i=1}^{N} \left(\frac{\tilde{a}(-\delta_i)
\tilde{a}(\delta_i)}{d'(\delta_i) d(-\delta_i)}\right) \]
The first order conditions (\ref{focV}) can be rewritten
as
\[\tilde{a}(\delta_i)=\frac{\tilde{a}(-\delta_i)^2} {e(\delta_i)},
 i=1,\ldots,N,~\tilde{a}\not=0.
\]
Substituting this and using $x_i=\tilde{a}(-\delta_i)$
it follows that
\[\| \frac{e(s)}{d(s)} - \frac{b(s)}{a(s)} \|_2^2 =
\sum_{i = 1}^N
\frac{x_i^3}{e(\delta_i)d'(\delta_i)d(-\delta_i)}.
\]
This shows 1.\\
Part 2 of the Lemma follows immediately from the fact that
the $L_2$ norm and the $H_2$ norm coincide for all elements in
$H_2.$ (See also the remark made above after the definition of
the $L_2-$norm).\\
Proof of part 3: Suppose that $a$ is not Hurwitz. Then it can be factored
uniquely as $a=a_1 a_2,$ where $a_1$ and $a_2$ are monic and
$a_1(s)$ and $a_2(-s)$ are Hurwitz polynomials in the variable
$s,$ with $\deg(a_1)<n.$ There are corresponding polynomials
$b_1, b_2$ with $\deg(b_1)<\deg(a_1)$ and $\deg(b_2)<\deg(a_2)$
such that
$\frac{b(s)}{a(s)}=\frac{b_1(s)}{a_1(s)}+\frac{b_2(s)}{a_2(s)}.$
As is well-known (and following from Cauchy's theorem in complex
analysis)
\[\frac{1}{2 \pi} \int_{-\infty}^{\infty} \frac{b_1(i\omega)
b_2(-i\omega)}{a_1(i\omega) a_2(-i\omega)} d\omega=0\]
and similarly
\[\frac{1}{2 \pi} \int_{-\infty}^{\infty} \frac{e(i\omega)
b_2(-i\omega)}{d(i\omega) a_2(-i\omega)} d\omega=0.\]
From this well-known orthogonality property
in $L_2$ it follows that
\[ \| \frac{e(s)}{d(s)} - \frac{b(s)}{a(s)} \|_2^2 =\]
\[ \| \frac{e(s)}{d(s)} - \frac{b_1(s)}{a_1(s)} \|_2^2 +
 \| \frac{b_2(s)}{a_2(s)} \|_2^2 \geq \]
\[ \| \frac{e(s)}{d(s)} - \frac{b_1(s)}{a_1(s)} \|_2^2.\]
This number is larger than the global minimum of
the function $f$ of Theorem \ref{th:H2critical}, because
 $\frac{b_1(s)}{a_1(s)}$ is the transfer function
of a system of order $< n.$ As noted before it is well-known
that the $H_2-$norm squared of the difference between the original
system and an approximant of order $<n,$ is always larger than
the global minimum of the $H_2-$norm squared of the difference
between the original
system and an approximant of order $n.$
This finishes the proof of part 3 and of the Theorem.\\
 \phantom{x} \hfill $\square$

For any complex polynomial $p \in {\bf C}[s]$ let $\bar{p}$ denote
the polynomial that is obtained when the coefficients
of $p$ are replaced by their complex conjugates. I.e.
$\bar{p}$ is the polynomial with the property that
$\bar{p}(r)=\overline{p(r)}$ for all $r \in {\bf R},$ where
$\bar{s}$ denotes the complex conjugate of a
complex number $s.$

\begin{lemma}
\label{complexsol}
%
Let $\frac{e}{d} \in \Sigma S_{N}$ have $N$ distinct poles
$\delta_1,\ldots,\delta_N \in {\bf C}.$

Let $a(s), b(s), q_0$ be a complex solution of  the polynomial
equations (\ref{bequation}),(\ref{focIII}). Then
$\bar{a}(s), \bar{b}(s), \overline{q_0}$ is also a solution.

The corresponding  numbers
$\sum_{i = 1}^N
\frac{q_0 a(-\delta_i)^3}{e(\delta_i)d'(\delta_i)d(-\delta_i)}$
and
$\sum_{i = 1}^N
\frac{\bar{q_0} \bar{a}(-\delta_i)^3}{e(\delta_i)d'(\delta_i)
d(-\delta_i)}$
form a complex conjugate pair.
In particular this implies that if one of these numbers is real
the numbers are equal.
\end{lemma}

{\bf Proof.}
Consider a complex solution $a(s), b(s),q_0$ of the first order
equations $e(s)a(s)-b(s)d(s)=a(-s)^2 q_0.$ Because polynomials are
 completely determined by their restriction to the real numbers,
an equivalent formulation of the first order equations is
$e(r)a(r)-b(r)d(r)=a(-r)^2 q_0$ for all $r \in {\bf R}.$
Complex conjugation of these equations gives
 $e(r)\bar{a}(r)-\bar{b}(r)d(r)=\bar{a}(-r)^2 \overline{q_0},$
which shows that $\bar{a}(s), \bar{b}(s), \overline{q_0}$ is also
a solution.

Because $h$ is a real polynomial with distinct zeros the set of
zeros of $h$ consists of an even number, $2l,$ say, of complex
solutions and $n-2l$ real solutions. The $2l$ complex solutions
can be partioned into $l$ pairs of complex conjugate solutions.
It is easy to see that for each real zero $\delta$ of $h,$
\[\frac{q_0 a(-\delta)^3}{e(\delta)d'(\delta)d(-\delta)}\]
and \[\frac{\bar{q_0} \bar{a}(-\delta)^3}{e(\delta)d'(\delta)
d(-\delta)}\] is a complex conjugate pair.
And if $\delta, \overline{\delta}$ is a complex conjugate pair
of zeros of $h,$ then the complex conjugate of
\[\frac{q_0 a(-\delta)^3}{e(\delta)d'(\delta)d(-\delta)}+
\frac{q_0 a(-\overline{\delta})^3}{e(\overline{\delta})
d'(\overline{\delta})d(-\overline{\delta})}\] is equal to
\[\frac{\overline{q_0} \bar{a}(-\overline{\delta})^3}{e(\overline{\delta})
d'(\overline{\delta})d(-\overline{\delta})}+
\frac{\overline{q_0} \bar{a}(-\delta)^3}{e(\delta)d'(\delta)d(-\delta)}.\]
Combining this it follows that
\[\sum_{i = 1}^N
\frac{q_0 a(-\delta_i)^3}{e(\delta_i)d'(\delta_i)d(-\delta_i)}\]
and
\[\sum_{i = 1}^N
\frac{\bar{q_0} \bar{a}(-\delta_i)^3}{e(\delta_i)d'(\delta_i)
d(-\delta_i)}\]
form a complex conjugate pair.
\\
 \phantom{x} \hfill $\square$

For ease of reference, let $\phi$ be defined by $\phi: \Xi \rightarrow
{\bf C}, x \mapsto \sum_{i = 1}^N
\frac{x_i^3}{e(\delta_i)d'(\delta_i)d(-\delta_i)}.$

Using the results above one can find the global minimum of the
criterion function as follows. For each of the at most $2^N-1$
elements of $\Xi,$ evaluate the numbers $\phi(x) \in {\bf C}.$ At
least one of these numbers will be real and positive. Let $k$
denote the number of distinct real positive numbers obtained in
this way and let us denote these numbers by $m_1,\ldots,m_k$ where
$m_1<\ldots<m_k.$ Consider the set  $\phi^{-1}(m_1).$ If each $\xi
\in \phi^{-1}(m_1)$ corresponds to a complex non-real solution
$a(s), b(s), q_0$ of  the polynomial equations
(\ref{bequation}),(\ref{focIII}), there must be an even number of
such solutions, as a result of Lemma (\ref{complexsol}). If any of
the solutions is real then according to Theorem \ref{realsol} the
global minimum is equal to $m_1$  and all real solutions $a(s),
b(s), q_0$ that correspond to this number are global approximants.
If none  of the solutions that correspond to $\xi \in
\phi^{-1}(m_1)$ are real then consider the set $\phi^{-1}(m_2).$
If any of the corresponding solutions $a(s), b(s), q_0$ is real
then $m_2$ is the global minimum, otherwise consider the solutions
 that correspond to  $m_3$ etc. One of the numbers
$m_1,\ldots,m_k$ is the global minimum and therefore the global
minimum will be found in this way. It follows from Theorem
\ref{realsol} that all real solutions $a(s), b(s), q_0$ that
correspond to the global minimum are in fact admissible, i.e.
$a(s)$ is Hurwitz and $a(s)$ and $b(s)$ are coprime.

{\em Remark.} Note that the function $\phi$ is a polynomial and
therefore continuous and smooth. Depending on the size of the
coefficients
  $\frac{1}{e(\delta_i)d'(\delta_i)d(-\delta_i)}$
a perturbation in $x$ due to numerical round-off error may
cause a limited perturbation in the corresponding value of $\phi.$
This implies that if the size of the coefficients just mentioned
is not too big, and the perturbations in $x$ are limited then the
effects of round-off error on the calculated critical values are
limited.
This can be contrasted with the possible effect of perturbations
on the calculation of the critical points. Especially if a
critical point $\frac{b(s)}{a(s)}\in \Sigma S_{N-1}$ has poles near the
imaginary axis, a small
perturbation may produce a denominator polynomial with one or more
right half-plane zeros, and therefore an inadmissible system,
outside the manifold $\Sigma S_{N-1}.$
Note that even if due to round-off error our algorithm would not produce
 a reliable global approximant, knowledge of the value of the global minimum
 of the criterion function could be used to evaluate the performance of
 other algorithms for the $H_2$ model order reduction problem.

{\em Remark.} The formula for $\phi$ in the Theorem can be used to
build the critical value matrix $A_F$ that was mentioned at the
end of the previous section, by taking the polynomial $f$
mentioned there equal to $\phi.$  Note that because $\phi$ is a
{\em polynomial}
 no matrix inversion is required in the calculation
of $A_F=\phi(A_{X_1},\ldots,A_{X_N}).$ The matrix $A_F$ can also
be built up by direct construction of the matrix of the
endomorphism $F$ with respect to the basis $B$ of monomials
defined earlier.

\section{Repeated poles}

In this section we briefly outline how the development is changed
if any of the poles of the original system are repeated, and
indicate the additional difficulty which arises in that case. For
simplicity of exposition we assume that one pole has multiplicity
two: $\delta_1=\delta_2$, and the other poles are distinct. In
this case \eqref{focV} gives only $N-1$ independent equations. An
additional equation is obtained by differentiating \eqref{focII},
which leads to
\begin{equation}\label{focrepeated}
    e(\delta_1)\tilde{a}'(\delta_1)+e'(\delta_1)\tilde{a}(\delta_1)
    = -2\tilde{a}(-\delta_1)\tilde{a}'(-\delta_1)
\end{equation}
(Note that we have used $d(\delta_1)=d'(\delta_1)=0$ here.) Taking
$x_1=\tilde{a}(-\delta_1)$, $x_2=\tilde{a}'(-\delta_1)$,
$x_i=\tilde{a}(-\delta_i)$ for $i=3,\ldots,N$, one obtains again a
system of $N$ quadratic polynomial equations in $x_1,\ldots,x_N$
representing the first-order conditions.

This system of equations will not yet be in Gr\"{o}bner basis
form, in contrast to the case of distinct poles. So at this point
it is necessary to employ Buchberger's algorithm to obtain a
Gr\"{o}bner basis for the corresponding ideal. Subsequently the
Stetter-M\"{o}ller matrix method can again be used to find the
critical points and hence the global optimum, provided that the
number of critical points is finite. As far as we are aware, there
is as yet no guarantee that this is the case.

If $\delta_1$ has multiplicity greater than two then higher-order
differentiation of \eqref{focII} is needed, but otherwise the
generalization is rather straightforward. If there are several
repeated poles a similar approach can be followed.

\section{Examples}
\label{sec:examples}

\subsection{General}
\label{sec:general}
This section presents two examples on solving the $H_2$ model
reduction problem and discusses several computational issues.

The following is an outline of the algorithm implemented:
\begin{enumerate}

\item For the given $N$-th order transfer function to be reduced,
construct the $N$-by-$N$ matrix $M$ (see equation (\ref{focVII})).

\item \label{alg:m2axi} For $i=1,...,N$, construct the $2^N$-by-$2^N$
  matrix $A_{X_i}$ from $M$ (see Theorem \ref{thm:fAXizero} and the
  following paragraph, and note
  that the reduction procedure of section \ref{quaddiag} is crucial in
  enabling this to be done).

\item \label{alg:eigvals}
  Compute the eigenvalues and eigenvectors of all the $A_{X_i}$'s.
  Assume, for simplicity, that each $A_{X_i}$ has a simple Jordan structure.
  Arrange these
  eigenvalues and eigenvectors such that the $j$-th eigenvector of
  $A_{X_{i_1}}$ corresponds to the $j$-th eigenvector of $A_{X_{i_2}}$
  for all $j = 1,...,2^N$ and $i_1,i_2 = 1,...,N$.  Letting $\xi_{i,j}$
  denote the $j$-th eigenvalue of $A_{X_i}$, form the $N$-tuples
  $(\xi_{1,j},\ldots,\xi_{N,j})$ for $j = 1,...,2^N$. Now each of these
  $N$-tuples contains the eigenvalues that correspond to one of the
  common eigenvectors of the set $\{A_{X_i}\}$. Our current implementation
  of this step uses numerical methods, so there are potential problems
  which can arise if eigenvalues and/or eigenvectors are repeated,
  or nearly so. We have not attempted to cope with all such eventualities.

\item \label{alg:apoly}
  Solve for $\tilde{a}_i$, using equation (\ref{atilantidel}), by making
  the association
  $$[\tilde{a}(-\delta_1),\ldots,\tilde{a}(-\delta_N)]=
  [\xi_{1,j},\ldots,\xi_{N,j}].$$
  Normalise the coefficients such that $a_{N-1} = 1$ to obtain $a_i$.
  Eliminate those polynomials $a(s)=s^{N-1}+a_{N-2}s^{N-2}+\ldots+a_0$
  which are not admissible pole polynomials of an approximating system,
  because they are not real Hurwitz.

\item \label{alg:ls} For each admissible pole polynomial $a(s)$, obtain the
  zero polynomial $b(s)$ from equation (\ref{bequation}). In practice
  the equation does not hold exactly, so a least-squares solution is found.

\end{enumerate}

All the above steps except that of computing eigenvalues and
eigenvectors can in principle be performed symbolically. Two different
implementations have been attempted and they differ only in whether
step \ref{alg:m2axi} is performed symbolically or numerically; note
that steps \ref{alg:eigvals} and \ref{alg:apoly} are done numerically
here.  For the symbolic implementation of step \ref{alg:m2axi}, the
$A_{X_i}$'s are computed from a symbolic definition of $M=[m_{jk}]$
using computer algebra software\footnote{In our case, Maple.} and the
resulting symbolic expressions for the $A_{X_i}$'s (see the Appendix)
are stored in a file to be read in by numerical software\footnote{In
  our case, Matlab.} later. This has the advantage that the symbolic
computation only has to be performed {\sl once} for a given model
order. Unfortunately, the length of these symbolic expressions soon
becomes very large as the model order increases; the size of the file
storing these expressions comes to more than 5 Mbytes for model order
equal to 7 and this thus presents a practical limit to this
implementation. Alternatively, due to the simplicity of the reduction
procedure (see section \ref{quaddiag}), step \ref{alg:m2axi} can be
implemented in a straightforward manner in a numerical
package${}^{\rm 2}$. In
this case, the highest model order that we could reduce is 9, which
involves storing 9 $512\times 512$ matrices, and we ran into memory
problems for model orders higher than this.
The computer we used was a Sun Ultra 10, 300 MHz
processor with 640 MByte RAM.

There are a number of numerical issues pertaining to this algorithm.
Some of these issues are well known, e.g.~possible ill-conditioning of
Vandermonde matrices and the computation of eigenvalues and
eigenvectors. These numerical problems will also cause difficulty in
later steps of the algorithm. For example, numerical error may cause
us to regard a real polynomial as complex in step \ref{alg:apoly} and
as a result, a true local minimum of the problem may be mistakenly
considered as inadmissible. The current
implementation of this algorithm does not strive to overcome nor
detect these problems. It is also beyond the scope of this paper to
give full numerical analysis of the proposed algorithm of this paper.
A rudimentary check that we have employed is to examine the
least-squares error in step \ref{alg:ls}; however, this error must be
interpreted with care as a small residual error does not necessarily
indicate an accurate solution \cite{GolubVL}.
Moreover, this check will not be able to tell us whether a correct
solution has been rejected. We have applied our algorithm to the three
third order systems that were investigated in \cite{HanMacAut96}
where a symbolic algorithm was used to reduce them to second order
systems.  In this case, symbolic computation ensures that all
stationary points of the problem are computed and we find that the
algorithm of this paper is able to find the same sets of critical
points as those reported in \cite{HanMacAut96}. This comparison may
indicate that our algorithm is likely to return the entire set of
stationary points when the model order is small.

\subsection{Example 1: An easily reduced system}
\label{sec:example1}
The system to be reduced is a 9th order transfer function and it
is the highest order model that we could reduce thus far. This system
has Hankel singular values $9,8,\ldots,2,1$ and its transfer function
is

\begin{math}
\frac{    8.4800 s^8  -2.5942 s^7 + 153.5350 s^6 + 38.8803 s^5 +
599.3205 s^4 + 196.3752 s^3 +  315.3021 s^2 +
 6.4558 s + 9.4478\times 10^{-5}}
    {s^9 +  2.1179 s^8 +   16.1278 s^7 +  25.6052 s^6 +  62.7884 s^5
      +  79.1895 s^4 +  42.6617 s^3 +
   32.5279 s^2 +  0.2514 s +  2.2495\times 10^{-6}}
\end{math}

The algorithm finds 8 admissible stationary points altogether. The best
approximant is

\begin{math}
\frac{8.4799 s^7  -2.5955  s^6 + 153.5327 s^5 +  38.8546 s^4 +
 599.3039 s^3 + 196.2798 s^2 + 315.2701 s +     6.4351   }
{s^8 + 2.1176 s^7 + 16.1275 s^6 + 25.6013 s^5 + 62.7850 s^4 +
79.1756 s^3 + 42.6527 s^2 +   32.5215 s +  0.2499}
\end{math}

and it gives $H_2$ model reduction error of 0.0344 and in comparison
with the norm of the original transfer function 8.8261, this gives a
relative error of 0.39\%. Note that the coefficients of this
approximant are very similar to those of the original transfer
function and this can be accounted for as follows: the original
transfer function has a pole at $-8.9582\times 10^{-6}$ and a zero at
$-1.4645\times 10^{-5}$. The model reduction algorithm appears to have
removed this very closely spaced pole-zero pair and to have left the
other poles and zeros nearly unchanged. The other seven approximants
give errors of 0.8703, 0.8707, 1.6463, 1.6466, 1.6536, 1.6538 and
1.6650.  Provided that all the stationary points of this optimisation
problem have been computed, then the solution that gives the minimum
error is in fact the global minimum of the problem. The other
stationary points may correspond to local minima, saddle points or
even local maxima.

\subsection{Example 2: A relaxation system}
\label{sec:example2}
The system to be reduced is taken from p.162 of \cite{ZhouDG} and is
given by
\begin{eqnarray}
G(s) & = & \sum_{j=1}^{N} \frac{\alpha^{2j}}{s+\alpha^{2j}} \; \;
\mbox{with } \alpha > 0.
\label{eq:relaxation}
\end{eqnarray}
It is shown in \cite{ZhouDG} that all the Hankel singular values of
this system tend to $\frac{1}{2}$ as $\alpha \rightarrow \infty$. On
the other hand, when $\alpha \approx 1$ and $N > 1$, the system is
close to non-minimality as $\alpha = 1$ gives rise to a first order
system. Our algorithm has numerical difficulty when $\alpha$ is chosen
either too large or too close to 1. In both cases, the Vandermonde
matrix becomes ill-conditioned: the rows contain entries of
drastically different magnitude in the first case and the poles are
too close to each other in the second.

Since the poles of this system are all real, this gives rise to a real
$M$ matrix and in turn real $A_{X_i}$'s. Due to the form of Gr\"obner
basis defined by $M$, zero is always an eigenvalue of $A_{X_i}$
(independent of whether $M$ is real or complex). Since the dimension
of $A_{X_i}$ is $2^N$ --- an even number --- and $A_{X_i}$ is real,
$A_{X_i}$ must have at least one other non-zero real eigenvalue. For
$\alpha$ close to zero or unity, we find in our examples there is a
real eigenvalue that is approximately zero and the eigenvectors
corresponding to this eigenvalue and the zero eigenvalue are almost
parallel to each other. This gives rise to difficulty in matching the
eigenvectors.

For model order $N = 5$, our algorithm succeeded in finding an
approximant for systems with $\alpha$ in the interval $[0.38,0.79]$
but failed in the intervals $(0,0.38)$ and $(0.79,1)$. For $\alpha$ in
the intervals $(0,0.38)$ and $(0.84,1)$, our algorithm returns no
solution as it either has difficulty in matching the eigenvectors or
has rejected the admissible solutions because they are not real
Hurwitz. Our algorithm does return a solution for $\alpha \in
(0.79,0.84]$ but a closer examination of the obtained approximant
shows that it is not a relaxation system. Since the system in equation
(\ref{eq:relaxation}) is a relaxation system and it is proved in
\cite{Baratchart93} that $H_2$ approximants of relaxation systems
are also relaxation systems, it implies that the solution given by our
algorithm for this range of $\alpha$ is unacceptable.

It is also shown in \cite{Baratchart93} that any stable relaxation
system, whose poles all have modulus smaller than
$\frac{1}{\sqrt{2}} \approx 0.707$, has only one admissible
solution of the first-order optimality conditions. For $\alpha =
0.78$, the largest pole is 0.6084 and there should therefore be
only one such solution. For this case our algorithm returns
precisely one admissible solution, in accordance with this theory.
It has absolute error 0.0334, which can be compared to the norm
1.6980 of the original system to give a relative error of 1.96\%.
The transfer function of this approximant is

\begin{math}
  \frac{1.4240 s^3 + 1.0946 s^2 + 0.2371 s + 0.0134} {s^4 + 1.1781 s^3
    + 0.4457 s^2 + 0.0627 s + 0.0028}.
\end{math}

which can be shown to be a relaxation system.

As an alternative to the algorithm described at the beginning of
this section, we have also treated Example 2 using an algorithm
based on building up the critical value matrix using
\eqref{thirddegree}. The same results were obtained with both
algorithms, except when $\alpha$ was very close to 1. For example
with $N=2$ and $\alpha=0.999$ the first algorithm continued to
give the correct result (which was checked using exact algebraic
calculation) but the second did not, because of numerical
imprecision.

\section{Conclusions}
The application of constructive algebra methods to the $H_2$ approximation
problem offers the possibility of guaranteed location of the globally
optimal approximant, despite the fact that this is a non-convex problem.
Furthermore, the location of this optimal approximant could, in principle,
be computed to any desired precision, by employing `symbolic' methods
throughout.

One can envision, however, that these methods could be used in
conjunction with more conventional numerical methods in a number
of ways, to obtain various precision/efficiency trade-offs. One
possibility is the one used by us to solve the examples in this
paper, namely to employ conventional numerical eigenvalue solvers
from the point at which the matrices $A_{X_i}$ have been
determined. Another possibility would be to use constructive
algebra methods to obtain an upper bound for the number of
admissible critical points, and/or the value of the criterion
function at the optimal approximant (which can be done without
computing the optimal approximant itself), and to use these
results to check the candidate optima obtained by a conventional
numerical optimization approach.

It should be kept in mind that constructive algebra also offers the
possibility of dealing with purely symbolic problem specifications ---
that is, of producing `generic' results (for all transfer functions of a
given order, say) rather than results for one specific system. This can be
done in principle, although in practice the complexity of the required
computations is well beyond current possibilities.

The commutative matrix approach which we have used to solve the
system of critical-point (polynomial) equations is currently the
subject of intense research activities in the computer algebra
community, and in the systems theory community \cite{BleHanPee03},
\cite{HanzonHazewinkel2006} with good prospects of much more
efficient algorithms being developed. We therefore expect that it
will soon be possible to approximate higher-order systems than the
ones we have been able to tackle in this paper, using essentially
the same methods. Also, we expect that such developments will make
constructive algebra methods attractive and feasible tools for
tackling a wider range of problems in systems and control theory.

\section*{Acknowledgements}
This research was supported by the British Council, the Dutch
Science Foundation NWO, and by the European Commission through the
European Research Network on System Identification ERNSI (Contract
ERB FMRX CT98 0206). Part of this research was done while B.Hanzon
was visiting Cambridge University Engineering Dept.

\end{document}